\documentclass{LMCS}

\def\doi{7 (4:02) 2011}
\lmcsheading%
{\doi}
{1--42}
{}
{}
{Nov.~\phantom28, 2010}
{Nov.~\phantom02, 2011}
{}   

\usepackage{enumerate,hyperref}
\usepackage{amsmath, amsthm, amssymb}
\usepackage{latexsym}
\usepackage{amsfonts}

\newcommand{\set}[2]{\ensuremath{\{#1 \ \mid \ #2\}}}

\newcommand{\om}{\ensuremath{\omega}}
\newcommand{\ep}{\ensuremath{\varepsilon}}

\newcommand{\iin}{\ensuremath{i \in \om }}
\newcommand{\n}{\ensuremath{n \in \om }}
\newcommand{\s}{\ensuremath{s \in \om }}

\newcommand{\ca}[1]{\ensuremath{ \mathcal{#1} }}

\newcommand{\eq}{\ensuremath{\Longleftrightarrow }}

\newcommand{\G}{\ensuremath{\Gamma }}

\newcommand{\empt}{\ensuremath{ \langle \cdot \rangle }}

\newcommand{\vg}[1]{}

\newcommand{\omseq}{\ensuremath{\om^{\star}}}

\begin{document}

\title[The complexity of the set of points of continuity]{The descriptive set-theoretic complexity of the set of points of continuity of a multi-valued function\rsuper*}

\author[V. Gregoriades]{Vassilios Gregoriades}
\address{Fachbereich Mathematik, Technische Universit\"{a}t Darmstadt, Schlossgartenstr. 7,\newline 64289 Darmstadt, Germany}
\email{gregoriades [at] mathematik
[dot] tu-darmstadt [dot] de}

\thanks{The author is currently a post-doctoral researcher at the Mathematics Department of TU Darmstadt
in the work-group of \textsc{Ulrich Kohlenbach} (TU Darmstadt),
whom the author would like to thank. The author owns special
thanks to \textsc{Martin Ziegler} (TU Darmstadt) for bringing to
the author's attention the problem which motivated this article
and for his valuable advice. Finally the author would like to thank all referees of this article (including the ones of the submission to the proceedings of the $7^{{\rm th}}$ International Conference on Computability and Complexity in Analysis) for their kind remarks, comments and suggestions.}

\titlecomment{{\lsuper*}An extended abstract of this article has been published in the proceedings of the $7^{{\rm th}}$ International Conference on Computability and Complexity in Analysis.}

\begin{abstract}
In this article we treat a notion of continuity for a multi-valued
function $F$ and we compute the descriptive set-theoretic
complexity of the set of all $x$ for which $F$ is continuous at
$x$. We give conditions under which the latter set is either a
$G_\delta$ set or the countable union of $G_\delta$ sets. Also we
provide a counterexample which shows that the latter result is
optimum under the same conditions. Moreover we prove that those
conditions are necessary in order to obtain that the set of points
of continuity of $F$ is Borel i.e., we show that if we drop some
of the previous conditions then there is a multi-valued function
$F$ whose graph is a Borel set and the set of points of continuity
of $F$ is not a Borel set. Finally we give some analogous results
regarding a stronger notion of continuity for a multi-valued
function. This article is motivated by a question of Martin Ziegler (TU Darmstadt).
\end{abstract}

\maketitle

\section{Introduction.}

\label{section introduction}

A \emph{multi-valued} function $F$ from a set $X$ to another set
$Y$ is any function from $X$ to the power set of $Y$ i.e., $F$
assigns sets to points. Such a function will be denoted by $F:X
\Rightarrow Y$. A multi-valued function $F: X \Rightarrow Y$ can
be identified with its \emph{graph} $Gr(F) \subseteq X \times Y$
which is defined by
\[
(x,y) \in Gr(F) \eq y \in F(x).
\]
This way we view $F$ as a subset of $X \times Y$. From now on we
assume that all given multi-valued functions are between metric
spaces and that they are \emph{total} i.e., if $F : X \Rightarrow
Y$ is given then $F(x) \neq \emptyset$ for all $x \in X$, in other
words the projection of $F$ along $Y$ is the whole space $X$.

There are various notions of continuity for multi-valued
functions, here we focus on two of those (see
\cite{brattka-hertling} Definition 2.1,
\cite{choquet_multi_valued_functions} pp. 70-71 and
\cite{adamovicz_multivalued_functions} p. 82, p. 93).

\begin{defi}

\normalfont

\label{definition of continuity of multivalued functions}

Let $(X,p)$ and $(Y,d)$ be metric spaces; a multi-valued function
$F: X \Rightarrow Y$ is \emph{continuous at $x$} if there is some
$y \in F(x)$ such that for all $\ep > 0$ there is some $\delta >0$
such that for all $x' \in B_p(x,\delta)$ there is some $y' \in
F(x')$ for which we have that $d(y,y') < \ep$.
\end{defi}

\begin{defi}

\normalfont

\label{definition of strong continuity of multivalued functions}

Let $(X,p)$ and $(Y,d)$ be metric spaces; a multi-valued function
$F: X \Rightarrow Y$ is \emph{strongly continuous at $x$} if for
all $y \in F(x)$ and for all $\ep > 0$ there is some $\delta
>0$ such that for all $x' \in B_p(x,\delta)$ there is some $y' \in
F(x')$ for which we have that $d(y,y') < \ep$.
\end{defi}

It is clear that both these notions generalize the classical
notion of continuity of functions. Moreover it is also clear that
the continuity/strong continuity of a multi-valued function is
preserved under distance functions which generate the same
topology.

The motivation of this article is the following question posed by
M. Ziegler in \cite{martin} (Question 63a)). It is well known that if we have a function $f: (X,p) \to (Y,d)$ then the set of points of continuity of $f$ is a $G_\delta$ subset of $X$; see for example 3.B in \cite{kechris}. So what can be said about the descriptive set-theoretic complexity of the set of points of continuity/strong continuity of a multi-valued function $F: (X,p) \Rightarrow (Y,d)$? In this article we present the answer for the case of continuity and present some analogous results for the case of strong continuity. The full answer for the latter case is still under investigation.

We proceed with basic terminology and notations. By \om \ we
denote the set of natural numbers (including the number $0$).
Suppose that $X$ and $Y$ are two topological spaces. We call a
function $f: X \to Y$ a \emph{topological isomorphism} between $X$
and $Y$ if the function $f$ is bijective, continuous and the
function $f^{-1}$ is continuous. We will also say that the space
$X$ is \emph{topologically isomorphic} with $Y$ if there exists a
topological isomorphism between $X$ and $Y$.

The \emph{Baire space} \ca{N} is the set of all sequences of
naturals i.e., $\ca{N} = \om^\om$ with the usual product topology.
We call the members of the Baire space \emph{fractions} and we
usually denote them by lower case Greek letters $\alpha, \beta$
etc. One choice of basic neighborhoods for the product topology on
\ca{N} is the collection of the following sets
\[
N(k_0,\dots,k_{n-1}) = \set{\alpha \in \ca{N}}{\alpha(0)=k_0,
\dots, \alpha(n-1) = k_{n-1}}
\]
where $k_0,\dots,k_{n-1} \in \om$. The set of ultimately constant
sequences  is clearly countable and dense in \ca{N}; thus the
latter is a separable space. For $\alpha, \beta \in \ca{N}$ with
$\alpha \neq \beta$ define
\[
d_{\ca{N}}(\alpha,\beta) = 1 / (\textrm{least} \ n \ [\alpha(n)
\neq \beta(n)] \ + 1).
\]
Also put $d_{\ca{N}}(\alpha,\alpha) =0$ for all $\alpha \in
\ca{N}$. It is not hard to see that the function $d_{\ca{N}}$ is a
complete distance function on \ca{N} which generates its topology.
From now on we think of the Baire space \ca{N} with this distance
function $d_{\ca{N}}$.

We denote by \ca{C} the subset of the Baire space \ca{N} which
consist of all sequences with values $0$ and $1$ i.e., $\ca{C} =
2^\om$. The set \ca{C} with the induced topology is a compact
space. It is not hard to see that \ca{C} is topologically
isomorphic with the usual Cantor set of the unit interval. This
result motivates us to call \ca{C} \emph{Cantor space}.

We denote by \omseq \ the set of all finite sequences of \om. If
$u \in \omseq$ then there are unique naturals
$n,k_0,\dots,k_{n-1}$ such that $u = (k_0,\dots,k_{n-1})$. The
\emph{length of $u$} is the previous natural $n$ and we denote it
by $lh(u)$. Also we write $u(i) = k_i$ for all $i < lh(u)$, so
that $u=(u(0),\dots,u(lh(u)-1))$. It is convenient to include the
\emph{empty sequence} in \omseq \ i.e., the one with zero length.
The latter will be denoted by \empt. So when we write
$u=(u(0),\dots,u(n-1))$ we will always mean in case where $n=0$
that $u = \empt$. If $u \in \omseq$ and $n \in \om$ we denote the
finite sequence $(u(0),\dots,u(lh(u)-1),n)$ by $u \ \hat{} \ (n)$.
We write $u \sqsubseteq v$ exactly when $lh(u) \leq lh(v)$ and
$u(i) = v(i)$ for all $i < lh(u)$ i.e., $u \sqsubseteq v$ means
that $v$ is an extension of $u$ or equivalently $u$ is an initial
segment of $v$.

A set $T \subseteq \omseq$ is called a \emph{tree} on \om \ if it
is closed under initial segments i.e.,
\[
v \in T \ \& \ u \sqsubseteq v \ \ \Longrightarrow \ \ u \in T.
\]
The members of a tree $T$ are called \emph{nodes} or
\emph{branches} of $T$. A tree $T$ is of \emph{finite branching}
if and only if for all $u \in T$ there are only finitely many $n
\in \om$ such that $u \ \hat{} \ (n) \in T$. A fraction $\alpha$
is an \emph{infinite branch of $T$} if and only if for all \n \ we
have that $(\alpha(0),\dots,\alpha(n-1)) \in T$. The \emph{body}
$[T]$ of a tree $T$ is the set of infinite branches of $T$.

For practical reasons when we refer to a tree $T$ we will always
assume that $T$ is not empty i.e., $\empt \in T$. Define
\[
Tr = \set{T \subseteq \omseq}{\textrm{the set $T$ is a tree on
\om}}.
\]
We may view every $T \in Tr$ as a member of $2^{\omseq}$ by
identifying $T$ with its characteristic function $\chi_T : \omseq
\to \{0,1\}$. Since the set \omseq \ is countable the space
$2^{\omseq}$ with the product topology is completely metrizable -
in fact it is topologically isomorphic with the Cantor space
\ca{C}. Moreover the set $Tr$ is a closed subset of \ca{C}. Indeed
let $T_i \in Tr$ for all \iin \ be such that $T_i \stackrel{i \to
\infty}{\longrightarrow}S$ for some $S \in 2^{\omseq}$; we will
prove that $S \in Tr$. From the hypothesis it follows that for all
$u \in \omseq$ there is some $i_0 \in \om$ such that for all $i
\geq i_0$ we have that
\[
u \in T_i \eq u \in S.
\]
Taking $u = \empt$ since $T_i \in Tr$ for all \iin \ we have that
$\empt \in S$ and so $S$ is not empty. Also if $u, v \in \omseq$
we find $i$ large enough so that $u \in T_i \eq u \in S$ and $v
\in T_i \eq v \in S$. So if $u \in S$ and $v \sqsubseteq u$ then
$u \in T_i$ and since $T_i$ is a tree we also have that $v \in
T_i$; hence $v \in S$. Therefore $S \in Tr$ and the set of trees
$Tr$ is closed in \ca{C}.

We make a final comment about trees. For any non-empty set $S$ of
finite sequences of naturals the tree $T$ which is \emph{generated
by $S$} is the following
\[
\set{u}{(\exists w \in S)[ u \sqsubseteq w]}
\]
i.e., the tree which is generated by $S$ is the tree which arises
by taking all initial segments of members of $S$.

Suppose that $X$ is a metric space. The family $\Sigma^0_1(X)$ is
the collection of all \emph{open} subsets of $X$. Inductively we
define the family $\Sigma^0_{n+1}(X)$ for $n \geq 1$as follows:
for $A \subseteq X$,
\[
A \in \Sigma^0_{n+1}(X) \ \eq \ A = \bigcup\limits_{\iin} A_i, \ \
\ \textrm{where $X \setminus A_i \in \Sigma^0_{k_i}(X)$ for some
$k_i \leq n$ for all \iin.}
\]
Put also
\[
\Pi^0_n(X) =\set{B \subseteq X}{X \setminus B \in \Sigma^0_n(X)}
\]
and $\Delta^0_n(X) = \Sigma^0_n(X) \cap \Pi^0_n(X)$ for all $n
\geq 1$. Notice that family $\Pi^0_1(X)$ is the collection of all
closed subsets of $X$, the family $\Sigma^0_2(X)$ is the
collection of all $F_\sigma$ subsets of $X$ and so on. By a simple
induction one can prove that $\Sigma^0_n(X) \cup \Pi^0_n(X)
\subseteq \Delta^0_{n+1}(X)$ for all $n \geq 1$. It is well known
that in case where $X$ admits a complete distance function and it
is an uncountable set then $\Sigma^0_n(X) \neq \Pi^0_n(X)$ for all
$n \geq 1$ and so the previous inclusion is a proper one for all
$n \geq 1$, (see \cite{kechris} and \cite{yiannis}).

The families $\Sigma^0_n(X), \Pi^0_n(X)$ are closed under finite
unions, finite intersections, and continuous pre-images i.e., if
$f: X \to Y$ is continuous and $B \subseteq Y$ is in
$\Sigma^0_n(Y)$ then $f^{-1}[B]$ is in $\Sigma^0_n(X)$. Moreover
it is clear that if $f: X \to Y$ is a topological isomorphism then
for all $A \subseteq X$ we have that $A \in \Sigma^0_n(X)$ if and
only if $f[A] \in \Sigma^0_n(Y)$ and similarly for $\Pi^0_n$
for all \n. Finally the family is $\Sigma^0_n(X)$ is closed under
countable unions, the family $\Pi^0_n(X)$ is closed under
countable intersections and the family $\Delta^0_n(X)$ is closed
under complements. We usually say that $A$ is in $\Sigma^0_n$ when
$X$ is easily understood from the context. It is clear that all
sets in $\Sigma^0_n$ are Borel sets.

We now deal with a bigger family of sets. Suppose that $X$ is
separable and that $X$ admits a complete distance function. A set
$A \subseteq X$ is in $\Sigma^1_1(X)$ or it is \emph{analytic} if
$A$ is the continuous image of a closed subset of a complete and
separable metric space.\footnote{The notion of an analytic set can
be treated in a more general context of spaces; however we prefer
to stay in the context of complete and separable metric spaces.}
It is well known that in the definition of analytic sets we may
replace the term ``continuous image" by ``Borel image" (i.e.,
image under a Borel measurable function) and/or the term ``closed
subset" by ``Borel set"", (see \cite{kechris} and \cite{yiannis}).
A set $B \subseteq X$ is in $\Pi^1_1(X)$ or it is
\emph{co-analytic} if the set $X \setminus B$ is analytic and $B$
is in $\Delta^1_1(X)$ or it is \emph{bi-analytic} if $B$ is both
analytic and co-analytic. It is well known that every Borel set is
analytic, hence every Borel set is bi-analytic. A classical
theorem states that a set $B$ is Borel if and only if it is
bi-analytic, (see \cite{yiannis} 2E.1 and 2E.2).

The families $\Sigma^1_1(X), \ \Pi^1_1(X)$ and $\Delta^1_1(X)$ are
closed under countable unions, countable intersections and
continuous pre-images. Moreover the family $\Sigma^1_1(X)$ is
closed under Borel images i.e., if $Y$ is a complete and separable
metric space, $f: X \to Y$ is a Borel measurable function and $A$
is an analytic subset of $X$ then $f[A]$ is an analytic subset of
$Y$. Finally if $X$ is uncountable we have that $\Sigma^1_1(X)
\neq \Pi^1_1(X)$ and in particular there is an analytic set which
is not Borel.

We can pursue this hierarchy further by defining the family $\Sigma^1_{n+1}(X)$ as the collection of all subsets of $X$ which are the continuous image of a $\Pi^1_n$ subset of a complete and separable metric space. Similarly one defines the family $\Pi^1_{n+1}(X)$ as the collection of all subsets of $X$ whose complement is in $\Sigma^1_{n+1}(X)$ and the family $\Delta^1_{n+1}(X)$ as the collection of all subsets of $X$ which belong both to $\Sigma^1_{n+1}(X)$ and $\Pi^1_{n+1}(X)$. By a simple
induction one can prove that $\Sigma^1_n(X) \cup \Pi^1_n(X)
\subseteq \Delta^1_{n+1}(X)$ for all $n \geq 1$. Also the analogous properties stated above are true. The reader may refer to \cite{yiannis} for more information on those classes.

The proofs of the forthcoming theorems make a substantial use of
techniques of Descriptive Set Theory which involve the use of many
quantifiers. Of course those quantifiers can be interpreted as
unions and intersections of sets and this is what we usually do in
order to prove that a given set is for example $\Pi^0_2$. There
are some cases though (for example in the proof of Theorem
\ref{theorem counterexample about S-3}) where this interpretation
becomes too complicated. In these cases it is better to think of a
given set $P$ as a relation in order to derive its complexity. The
reader can consult section 1C of \cite{yiannis} on how one can
make computations with relations.

\section{Results about the set of points of continuity of a multi-valued function.}

\label{Results about the set of points of continuity}

In this section we examine the set of points of continuity (as given in Definition \ref{definition of continuity of multivalued functions}) of a multi-valued function. We can classify our results by two groups: the ``low-level" group where the set of points of continuity is at most a $\Sigma^0_3$ set and the ``higher level group" where the set of points of continuity is not even a Borel set, while on the other hand under some reasonable assumptions it is an analytic set.

\subsection*{The case of low-level classification in the Borel hierarchy.}

We begin with some positive results regarding the set of points of
continuity of a multi-valued function $F$. Recall that a
topological space $Y$ is \emph{exhaustible by compact sets} if
there is a sequence $(K_n)_{\n}$ of compact subsets of $Y$ such
that every $K_n$ is contained in the interior of $K_{n+1}$ and $Y
= \bigcup\limits_{\n} K_n$. Notice the lack of any hypothesis
about the set $F$ in the next theorem.

\begin{thm}

\label{theorem G-d when compact and Sigma-3 when sigma-compact}

Let $(X,p)$ and $(Y,d)$ be metric spaces with $(Y,d)$ being
separable and let $F: X \Rightarrow Y$ be a multi-valued function.

\begin{enumerate}[\em(a)]

\item If  the set $F(x)$ is compact for all $x \in X$ then the set of
points of continuity of $F$ is $\Pi^0_2$ i.e., $G_\delta$.
\item If $Y$ is exhaustible by compact sets and
the set $F(x)$ is closed for all $x \in X$, then the set of points
of continuity of $F$ is $\Sigma^0_3$.
\end{enumerate}

\end{thm}

\proof
For $y \in Y$ and for a non-empty $A \subseteq Y$ we denote by
$d(y,A)$ the non-negative number $\inf\set{d(y,z)}{z \in A}$. It
is clear that when $A$ is compact then $d(y,A) = d(y,z_0)$ for
some $z_0 \in A$. It is useful to adopt the notation
$d(y,\emptyset) = 1$. We denote by $P$ the set of $x$'s for which
$F$ is continuous at $x$ and also we fix a sequence $(y_s)_{s \in
\om}$ in $Y$ for which the set $\set{y_s}{\s}$ is dense in $Y$.

We begin with assertion (a). We claim that
\[
x \in P \eq (\forall n)(\exists s) \
\inf\set{\sup\set{d(y_s,F(x'))}{x'\in B_p(x,\delta)}}{\delta
> 0} < \frac{1}{n+1} \ \ \ \ \ \ (*)
\]
Assume that $x \in P$ and let $\n$. There exists $y \in F(x)$ such
that for all $\ep > 0$ there is some $\delta >0$ such that for all
$x' \in B_p(x,\delta)$ there is some $y' \in F(x')$ for which we
have that $d(y,y') < \ep$. Choose $s \in \om$ such that
\[
d(y,y_s) < \frac{1}{2(n+1)}.
\]
We will show that $\inf\set{\sup\set{d(y_s,F(x'))}{x'\in
B_p(x,\delta)}}{\delta > 0} < \frac{1}{n+1}$. For this it is
enough to find some $\delta > 0$ for which we have that
\[
\sup\set{d(y_s,F(x'))}{x'\in B_p(x,\delta)} < \frac{1}{n+1} \ \ \
\ \ (1)
\]
Let $\ep = \frac{1}{2(n+1)} - d(y,y_s) > 0$ and choose $\delta >
0$ such that for all $x' \in B_p(x,\delta)$ there is some $y' \in
F(x')$ for which we have that $d(y,y') < \ep$. We claim that this
$\delta$ satisfies $(1)$. Indeed for $x' \in B_p(x,\delta)$ choose
a $y' \in F(x')$ which satisfies $d(y,y') < \ep$ and compute
\begin{eqnarray*}
d(y_s,F(x')) \leq d(y_s,y') \leq d(y_s,y) + d(y,y')< d(y_s,y) + \ep = \frac{1}{2(n+1)}.
\end{eqnarray*}
So $\sup\set{d(y_s,F(x'))}{x' \in B_p(x,\delta)} \leq
\frac{1}{2(n+1)} < \frac{1}{n+1}$ and we have proved the
left-to-right-hand direction.

Now we deal with the inverse direction. Assume the right-hand
condition in $(*)$. It follows that there exists a sequence
$(s_n)_{\n}$ of naturals such that
\[
\inf\set{\sup\set{d(y_{s_n},F(x'))}{x' \in B_p(x,\delta)}}{\delta
> 0} < \frac{1}{n+1} \ \ \ \ \ \ \ (2)
\]
for all \n. So for all \n \ there is some $\delta_n > 0$ such that
\[
\sup\set{d(y_{s_n},F(x'))}{x' \in B_p(x,\delta_n)} <
\frac{1}{n+1}.
\]
In particular we have that $d(y_{s_n},F(x)) < \frac{1}{n+1}$ for
all \n. There is some $z_n \in F(x)$ such
that $d(y_{s_n},z_n) < \frac{1}{n+1}$, for all
\n. From the compactness of $F(x)$ there is a subsequence
$(z_{k_n})_{\n}$ of $(z_n)_{\n}$ which converges to some $y \in
F(x)$. And since $d(y_{s_n},z_n) < \frac{1}{n+1}$ for all \n \ we
also have that the sequence $(y_{s_{k_n}})_{\n}$ is convergent to
$y$.

We claim that this $y \in F(x)$ witnesses that $x \in P$. Suppose
that we are given $\ep > 0$. Since $y_{s_{k_n}}
\stackrel{\n}{\longrightarrow} y$ there is some $n_0$ such that
for all $n \geq n_0$ we have that
\[
d(y_{s_{k_n}},y) < \frac{\ep}{2}.
\]
Choose some $n_1 \geq n_0$ for which we have that $\frac{1}{n_1} <
\frac{\ep}{2}$. Applying $(2)$ to the natural $k_{n_1}$ we obtain
that there is some $\delta > 0$ such that
\[
\sup\set{d(y_{s_{k_{n_1}}},F(x'))}{x' \in B_p(x,\delta)} <
\frac{1}{k_{n_1}+1}.
\]
Let $x' \in B_p(x,\delta)$; since $F(x')$ is compact there is some
$y' \in F(x')$ such that $d(y_{s_{k_{n_1}}},F(x')) =
d(y_{s_{k_{n_1}}},y')$. We will prove that $d(y,y') < \ep$. Indeed
\begin{eqnarray*}
d(y,y') &\leq& d(y,y_{s_{k_{n_1}}}) + d(y_{s_{k_{n_1}}},y') \\
&<& \frac{\ep}{2} + d(y_{s_{k_{n_1}}},F(x')) \\
&<& \frac{\ep}{2} + \frac{1}{k_{n_1}+1} < \frac{\ep}{2} + \frac{1}{n_1}  < \ep.
\end{eqnarray*}
Thus we have proved $(*)$. Now for fixed $n, s \in \om$ we define
the set $Q_{n,s}$ as follows:
\[
x \in Q_{n,s} \eq \inf\set{\sup\set{d(y_s,F(x'))}{x'\in
B_p(x,\delta)}}{\delta > 0} < \frac{1}{n+1},
\]
so that from $(*)$ we have that $P = \bigcap\limits_{n}
\bigcup\limits_{s} Q_{n,s}$. The proof will be complete as long as
we prove that the set $Q_{n,s}$ is open for all $n, s$. The latter
assertion is immediate from the fact that the ball $B_p(x,\delta)$
is open. To see this fix $n, s \in \om$ and let $x \in Q_{n,s}$.
Then there is some $\delta > 0$ such that
\[
\sup\set{d(y_s,F(x'))}{x'\in B_p(x,\delta)} < \frac{1}{n+1}.
\]
We claim that $B_p(x,\delta) \subseteq Q_{n,s}$. Indeed let $w \in
B_p(x,\delta)$; since the latter set is open there is some
$\delta_w > 0$ such that $B_p(w,\delta_w) \subseteq
B_p(x,\delta)$. Hence
\[
\set{d(y_s,F(x'))}{x'\in B_p(w,\delta_w)} \subseteq
\set{d(y_s,F(x'))}{x'\in B_p(x,\delta)}
\]
and so
\[
\sup\set{d(y_s,F(x'))}{x'\in B_p(w,\delta_w)} \leq
\sup\set{d(y_s,F(x'))}{x'\in B_p(x,\delta)} < \frac{1}{n+1}.
\]
Therefore
\begin{eqnarray*}
&& \inf\set{\sup\set{d(y_s,F(x'))}{x'\in B_p(w,\delta)}}{\delta > 0} \\
&\leq& \sup\set{d(y_s,F(x'))}{x'\in B_p(w,\delta_w)}\\
&<& \frac{1}{n+1}
\end{eqnarray*}
i.e., $w \in Q_{n,s}$ and we are done.

\mbox{}

The following remark is motivated by the comments of one of the referees of this article, who noticed that we can eliminate one instance of our hypothesis about compactness. The author would like to express his thanks.

\textbf{Remark.} The reader may notice that the only place in the whole proof of (a) where the hypothesis
about compactness is needed is in the inverse direction of $(*)$. The use of compactness of $F(x')$ for $x' \in B_p(x,\delta)$ can be avoided by choosing $y' \in F(x')$ such that $d(y_{s_{k_{n_1}}},y') < d(y_{s_{k_{n_1}}},F(x')) + \frac{\ep}{3}$ and by replacing the previous occurrences of $\frac{\ep}{2}$ by $\frac{\ep}{3}$. In conclusion what is proved here is the slightly stronger result: if $Y$ is separable, $x \in X$ and $F(x)$ is compact, then $F$ is continuous at $x$ exactly when the right hand of $(*)$ holds.

\mbox{}

We now prove (b). Let $(K_m)_{m \in \om}$ be a sequence of compact
subsets of $Y$ such that $K_m \subseteq K_{m+1}^{\circ}$ for all
$m \in \om$ and $Y = \bigcup\limits_{m \in \om} K_m$. We claim the
following modified version of $(*)$:
\begin{eqnarray*}
x \in P &\eq& (\exists m)(\forall n)(\exists s) \\ &&
\inf\set{\sup\set{d(y_s,F(x')\cap K_m)}{x'\in
B_p(x,\delta)}}{\delta
> 0 } < \frac{1}{n+1} \ \ \ \ \ \ (\dagger)
\end{eqnarray*}
The proof of $(\dagger)$ is similar to the proof of $(*)$; however
some mild changes are needed. We give a sketch of this proof. For
the left-to-right-hand direction let $y \in Y$ which witnesses
that $x$ is a point of continuity of $F$. Let $m$ be such that $y
\in K_m^{\circ}$ and choose $r > 0$ for which we have that
$B_d(y,r) \subseteq K_m$. Proceed as previously but instead of
taking $\ep = \frac{1}{2(n+1)} - d(y,y_s)$ take $\ep =
\min\{\frac{1}{2(n+1)} - d(y,y_s), r\} > 0$. Now notice that for
any $y' \in Y$ which satisfies $d(y,y') < \ep$ we have that $y'
\in K_m$. So the resulting $y' \in F(x')$ which witnesses that
$d(y_s,F(x')) \leq \frac{1}{2(n+1)}$ is also a member of $K_m$. It
follows that $d(y_s,F(x') \cap K_m) \leq \frac{1}{2(n+1)}$ as
well. For the inverse direction replace $F(x)$ and $F(x')$ with
the compact sets $F(x) \cap K_m$ and $F(x') \cap K_m$
respectively. Notice that when $d(y_s,F(x') \cap K_m) < 1$ then
-by our convention $d(y,\emptyset) =  1$- we have that $F(x') \cap
K_m \neq \emptyset$. Having proved $(\dagger)$ one defines the set
$R_{m,n,s}$ ($m,n,s \in \om$) as follows:
\[
x \in R_{m,n,s} \eq \inf\set{\sup\set{d(y_s,F(x') \cap K_m)}{x'\in
B_p(x,\delta)}}{\delta > 0} < \frac{1}{n+1}
\]
and as before it is easy to see that the set $R_{m,n,s}$ is open.
From $(\dagger)$ we have that $P =
\bigcup\limits_{m}\bigcap\limits_{n}\bigcup\limits_{s}R_{m,n,s}$
and so the set $P$ is $\Sigma^0_3$.
\qed

Readers who find that the convention $d(y,\emptyset) = 1$ in the previous proof is rather artificial may notice that the condition
\[
\inf\set{\sup\set{d(y_s,F(x')\cap K_m)}{x'\in
B_p(x,\delta)}}{\delta
> 0} < \frac{1}{n+1}
\]
simply means that
\[
(\exists \delta > 0)(\exists q \in \mathbb{Q}^+)(\forall x' \in
B_p(x,\delta))[ \ F(x') \cap K_m \neq \emptyset \ \& \ d(y_s,F(x')
\cap K_m) < q < \frac{1}{n+1} \ ].
\]
Theorem \ref{theorem G-d when compact and Sigma-3 when
sigma-compact} has an interesting corollary which answers to Question
63a) posed by M. Ziegler in \cite{martin}.

\begin{cor}

\label{corollary from theorem G-d and Sigma-3}

Suppose that $X$ is a metric space and that $F: X \Rightarrow
\mathbb{R}^m$ is a multi-valued function such that the set $F(x)$
is closed for all $x \in X$.

\begin{enumerate}[\em(a)]

\item The set of the points of continuity of $F$ is $\Sigma^0_3$.

\item If moreover the set $F(x)$ is bounded for all $x \in X$ then the
set of points of continuity of $F$ is $\Pi^0_2$.
\qed
\end{enumerate}
\end{cor}

\noindent
We now show that the results of Theorem \ref{theorem G-d when
compact and Sigma-3 when sigma-compact} are optimum. It is well
known that there are functions $f: [0,1] \to \mathbb{R}$ for which
the set of points of continuity is not $F_\sigma$. Therefore the
$\Pi^0_2$-answer is the best one can get. Thus we only need to
deal with the $\Sigma^0_3$-answer. The following lemmas, although
being straightforward from the definitions, will prove an elegant
tool for the constructions that will follow.

\begin{lem}

\label{lemma about extension}

Suppose that $(X_0,p_0)$, $(X_1,p_1)$, $(Y,d)$ are metric spaces
and that there exists a function $f: X_0 \to X_1$ such that $f[X_0]$ is closed and $f: X_0 \to f[X_0]$ is a topological isomorphism. Assume
that we are given a multi-valued function $F: X_0 \Rightarrow Y$.
Define the multi-valued function $\tilde{F}: X_1 \Rightarrow Y$ as
follows:
\[
\tilde{F}(x_1) = F(x_0) \ \ \ \textrm{if $x_1 = f(x_0)$ for some
$x_0 \in X_0$ and} \ \ \ \tilde{F}(x_1) = Y \ \ \
\textrm{otherwise}.
\]
Then
\begin{enumerate}[\em(1)]

\item $\tilde{F}$ is continuous at $x_1$ if and only if either $x_1 \not
\in f[X_0]$ or $x_1=f(x_0)$ and $F$ is continuous at $x_0$. Hence
if we denote by $P_0$ and $P_1$ the set of points of continuity of
$F$ and $\tilde{F}$ respectively we have that
\[
P_1 = f[P_0] \cup (X_1 \setminus f[X_0]).
\]
\item If $\Gamma$ is any of the classes $\Sigma^0_n, \Pi^0_n$, with
$n\geq2$ or $\Delta^1_1$ then
\[
P_1 \in \Gamma \eq P_0 \in \Gamma.
\]
In particular if the set of points of continuity of $F$ is not
$\Pi^0_3$ (Borel) then the set of points of continuity of
$\tilde{F}$ is not $\Pi^0_3$ (Borel respectively).

Moreover the sets $F$ and $\tilde{F}$ as subsets of $X_0 \times Y$ and $X_1 \times Y$ respectively satisfy
\[
F \in \Gamma \eq \tilde{F} \in \Gamma.
\]

\item If $F(x_0)$ is a closed subset of $Y$ for all $x_0 \in X_0$
then $\tilde{F}(x_1)$ is also a closed subset of $Y$ for all $x_1
\in X_1$.
\end{enumerate}

\end{lem}

\proof
Clearly we only need to prove the first assertion. Suppose that
$x_1 \in X_1 \setminus f[X_0]$. Since the set $f[X_0]$ is closed
then there is some $\delta > 0$ such that $B_{p_1}(x_1,\delta)
\cap f[X_0] = \emptyset$. So for all $x_1' \in
B_{p_1}(x_1,\delta)$ we have that $\tilde{F}(x_1') =
\tilde{F}(x_1) = Y$ and hence $\tilde{F}$ is continuous at $x_1$.

Assume that $x_1 = f(x_0)$ where $x_0$ is a point of continuity of
$F$. Let $y \in Y$ which witnesses the continuity of $F$ at $x_0$.
We will show that this $y$ witnesses the continuity of $\tilde{F}$
at $x_1$. First notice that $F(x_0) = \tilde{F}(x_1)$ and so $y
\in \tilde{F}(x_1)$. Assume now that we are given $\ep > 0$.
Choose $\delta_0 > 0$ such that for all $x_0' \in
B_{p_0}(x_0,\delta_0)$ there is some $y' \in F(x_0')$ such that
$d(y,y') < \ep$. Since the function $f^{-1}: f[X_0] \to X_0$ is
continuous for this $\delta_0 > 0$ there is some $\delta_1 > 0$
such that for all $x_1' \in B_{p_1}(x_1,\delta) \cap f[X_0]$ we
have that $f^{-1}(x_1) \in B_{p_0}(x_0,\delta_0)$.

Let $x_1' \in B_{p_1}(x_1,\delta_1)$. We claim that there is some
$y' \in \tilde{F}(x_1')$ such that $d(y,y') < \ep$. We take cases:
if $x_1' \not \in f[X_0]$ then $\tilde{F}(x_1') = Y$ and so we may
take $y' = y$; if $x_1' = f(x_0')$ for some $x_0'\in X_0$, then
$x_0'$ belongs to $B_{p_0}(x_0,\delta_0)$ and so there is some $y'
\in F(x_0') = \tilde{F}(x_1')$ such that $d(y,y') < \ep$.

Now assume that $x_1 = f(x_0)$ and $x_0$ is not a point of
continuity of $F$. Let $y \in \tilde{F}(x_1) = F(x_0)$. Since $F$
is not continuous at $x_0$ there is some $\ep^{*} > 0$ such that
for all $\delta_0 > 0$ there is some $x_0' \in
B_{p_0}(x_0,\delta_0)$ such that for all $y' \in F(x_0')$ we have
that $d(y,y') \geq \ep^{*}$. We claim that this $\ep^{*}$
witnesses that $\tilde{F}$ is not continuous at $x_1$.

Let $\delta_1$; since $f$ is continuous there is some $\delta_0 >
0$ such that for all $x_0' \in B_{p_0}(x_0,\delta_0)$ we have that
$f(x_0') \in B_{p_1}(x_1,\delta_1)$. Choose $x_0' \in
B_{p_0}(x_0,\delta_0)$ such that for all $y' \in F(x_0')$ we have
that $d(y,y') \geq \ep^{*}$. Let $x_1' = f(x_0')$; then $x_1' \in
B_{p_1}(x_1,\delta_1)$ and also for all $y' \in \tilde{F}(x_1') =
F(x_0')$ we have that $d(y,y') \geq \ep^{*}$.
\qed

Lemma \ref{lemma about extension} has a cute corollary which might
be regarded as the multi-valued analogue of the Tietze Extension
Theorem.

\begin{cor}

\label{corollary analogue of Tietze extensionality}

Every continuous multi-valued function which is defined on a
closed subset of a metric space can be extended continuously on
the whole space.\qed
\end{cor}

\begin{lem}

\label{lemma about composition}

Let $X, Y, Z$ be metric spaces, $F: X  \Rightarrow Y$ be a multi-valued function and $\pi: Y \to Z$ be a topological isomorphism between $Y$ and $\pi[Y]$. Define the composition $\pi \circ F : X \Rightarrow Z:$
\[
(\pi \circ F)(x) = \pi[F(x)], \ \ x \in X.
\]
The following hold.
\begin{enumerate}[\em(1)]

\item A point $x \in X$ is a point of continuity of $F$ if and only if $x$ is a point of continuity of $\pi \circ F$;

\item If the set $\pi[Y]$ is closed and the set $F(x)$ is closed for some $x \in X$ then the set $(\pi \circ F)(x)$ is also closed.

\item If $F$ is a Borel subset of $X \times Y$ then $\pi \circ F$ is a Borel subset of $X \times Z$.
\end{enumerate}

\end{lem}

\proof
The first assertion is proved as Lemma \ref{lemma about extension} using the following remarks. For all $z = \pi(y)$ and for all $\ep^* > 0$ there is some $\ep > 0$ such that for all $y' \in B_Y(y,\ep)$ we have that $\pi(y') \in B_Z(z,\ep^*)$; and for all $y \in Y$ and for all $\ep > 0$ there is some $\ep^* > 0$ such that for all $z' \in B_Z(\pi(y),\ep^*) \cap \pi[Y]$ we have that $\pi^{-1}(z) \in B_Y(y,\ep)$.

For the second assertion we know that if $F(x)$ is closed then since $\pi$ is a topological isomorphism the set $\pi[F(x)]$ is closed in $\pi[Y]$. So if moreover $\pi[Y]$ is closed in $Z$ it follows that $\pi[F(x)]$ is closed in $Z$.

For the third assertion notice that
\[
(x,z) \in \pi \circ F \eq z \in \pi[Y] \ \& \ (x,\pi^{-1}(z)) \in F.
\]
and so $\pi \circ F = (X \times \pi[Y]) \cap h^{-1}[F]$, where $h(x,y)=(x,\pi(y))$. From 2E.9 of \cite{yiannis} we have that the set $\pi[Y]$ is Borel and hence $\pi \circ F$ is Borel.
\qed

\begin{thm}

\label{theorem counterexample about S-3}

There is a multi-valued function $F: [0,1] \to \mathbb{R}$ such
that the set $F(x)$ is closed for all $x$, the set $F$ is a
$\Pi^0_2$ subset of $[0,1] \times \mathbb{R}$ and the set of
points of continuity of $F$ is not $\Pi^0_3$. Therefore the
$\Sigma^0_3$-answer is the best possible for a multi-valued
function $F$ from $[0,1]$ to $\mathbb{R}$ even if $F$ is below the
$\Sigma^0_3$-level.
\end{thm}

Using the fact that the set of points of continuity of a function $f: (X,p) \to (Y,d)$ is a $G_\delta$ set, we conclude to the following.

\begin{cor}

\label{corollary continuity is non-metrizable}

The notion of continuity which is given in Definition \ref{definition of continuity of multivalued functions} for a multi-valued function $F: (X,p) \Rightarrow (Y,d)$ cannot be reduced to the usual notion of continuity if we view $F$ as a single-valued function from $X$ to the powerset $\mathcal{P}(Y)$ of $Y$. To be more specific \emph{it is not true in general} that if we are given a multi-valued function $F: (X,p) \Rightarrow (Y,d)$ there is a metric $d_F$ on the set $F[X] = \set{F(x)}{x \in X}$ such that for all $x \in X$, the function $F: (X,p) \to (F[X],d_F)$ is continuous at $x$ in the usual sense exactly when $F$ is continuous at $x$ in the sense of Definition \ref{definition of continuity of multivalued functions}.
\qed
\end{cor}

Now we proceed with the proof of Theorem \ref{theorem counterexample about S-3}.

\proof
Since $2^{\om \times \om}$ is isomorphic with the space \ca{C},
from Lemma \ref{lemma about extension} it is enough to define a
multi-valued function $F: 2^{\om \times \om} \Rightarrow
\mathbb{R}$ such that the set $F(\gamma)$ is closed for all
$\gamma \in 2^{\om \times \om}$, the set $F$ is $\Pi^0_2$ and also the set of points of continuity of $F$ is not $\Pi^0_3$.

A typical example of a $\Sigma^0_3$ set which is not $\Pi^0_3$ is
the following:
\[
R = \set{\gamma \in 2^{\om \times \om}}{(\exists m)(\forall
n)(\exists s \geq n)[ \gamma(m,s) = 1]},
\]
(see \cite{kechris} 23.1). For all $m \in  \om$ define $R_m$ as
the set of all $\gamma \in 2^{\om \times \om}$ such that for all
$n$ there is some $s \geq n$ such that $\gamma(m,s) = 1$, so that
$R = \bigcup\limits_{m} R_m$.

If for some $m$ and $\gamma$ we have that $\gamma$ is not a member of $R_m$, then there is some $n$ such that for all $s \geq n$ we have
that $\gamma(m,s) = 0$. For $\gamma \not \in R_m$ put
\[
n(\gamma,m) = \ \textit{the least \n} \ [\textit{for all $s \geq n$ we have
that} \ \gamma(m,s) = 0] + 1.
\]
Define $F: 2^{\om \times \om} \Rightarrow \mathbb{R}$ as follows
\[
F(\gamma) = \set{m + \frac{1}{n(\gamma,m)+1}}{\gamma \not \in R_m}
\cup \set{m}{\gamma \in R_m}.
\]
Notice that $n(\gamma,m) + 1\geq 2$
for all $\gamma \not \in R_m$. Clearly the set $F(\gamma)$ is
closed for all $\gamma \in 2^{\om \times \om}$. First we prove
that the set $F$ is a $\Pi^0_2$ subset of $2^{\om \times \om}
\times \mathbb{R}$. Indeed
\begin{eqnarray*}
(\gamma,y) \in F &\eq& (\exists m \leq y) \{ \ \textrm{either} \ [ \ y =m \ \& \ (\forall
n)(\exists s \geq n)[ \gamma(m,s) = 1] \ ]\\
&& \textrm{or} \ [ \ y \neq m \ \& \ \frac{1}{y-m}-2 \in \om \ \& \ (\forall s \geq \frac{1}{y-m}-2)[ \gamma(m,s) = 0] \\
&&  \& \ (\forall i < \frac{1}{y-m}-2)(\exists s > i)[\gamma(m,s) = 1] \ ] \ \}.
\end{eqnarray*}

We now
claim that
\[
\gamma \in R \ \eq \ F \ \textrm{is continuous at $\gamma$}.
\]
Let $\gamma \in R$. Fix some $m \in \om$ such that for all $n$
there is some $s \geq n$ such that $\gamma(m,s) = 1$ i.e., $\gamma \in R_m$. We claim that the natural $m$ which belongs to $F(\gamma)$ witnesses
the continuity of $F$ at $\gamma$. Let $\ep > 0$ and $n \in \om$
such that $\frac{1}{n+1} < \ep$. Since $\gamma \in R_m$ for this $n$ we can choose some $s_n \geq n$ such that $\gamma(m,s_n) = 1$. Define
\[
W = \set{\beta \in 2^{\om \times \om}}{\beta(m,s_n) =
\gamma(m,s_n)=1}.
\]
Clearly $W$ is a basic neighborhood of $\gamma$. We will prove
that for all $\beta \in W$ there is some $y' \in F(\beta)$ such
that $|m- y'| < \ep$. Indeed if $\beta \in W$ and $\beta \in R_m$
then $m \in F(\beta)$ so one can take $y'= m$; if $\beta \not \in
R_m$ since $\beta(m,s_n) = \gamma(m,s_n) = 1$ we have that
$n(\beta,m) > s_n$. Hence for $y' = m + \frac{1}{n(\beta,m) + 1}$
we have that $y' \in F(\beta)$ (because $\beta \not \in R_m$) and
also that $|m - y'| = \frac{1}{n(\beta,m) + 1} < \frac{1}{s_n + 1}
\leq \frac{1}{n + 1} < \ep$. Therefore $F$ is continuous at
$\gamma$.

Assume now that $\gamma \not \in R$; clearly $\gamma \not \in R_m$
for all $m \in \om$. From the definition it follows that $F(\gamma) = \set{m +
\frac{1}{n(\gamma,m)+1}}{m \in \om}$. Let any $y \in F(\gamma)$;
then $y = m + \frac{1}{n(\gamma,m) + 1}$ for some $m \in \om$.
Since $n(\gamma,m) + 1\geq 2$ it follows that $y = m +
\frac{1}{n(\gamma,m) + 1} \leq m + \frac{1}{2}$. Hence for $\ep =
\frac{1}{2} \cdot \frac{1}{n(\gamma,m) + 1}> 0$ we have that $|y -
y'| \geq \ep$ for all $y' \in \om$. We will show that for all
basic neighborhoods $W$ of $\gamma$ there is some $\beta \in W$
such that for all $y' \in F(\beta)$ we have that $|y - y'| \geq
\ep$. Let $V$ be any basic neighborhood of $\gamma$. Then there
are some naturals $m_0, \dots, m_k$, $s_0, \dots, s_k$ such that
\[
V = \set{\beta \in 2^{\om \times \om}}{\forall i \leq k \
\beta(m_i,s_i) = \gamma(m_i,s_i)}.
\]
Define $\beta \in 2^{\om \times \om}$ as follows:
$\beta(m_i,s_i) = \gamma(m_i,s_i)$ for all $i \leq k$; and
$\beta(m,s) = 1$ in any other case. It is clear that $\beta \in V$ and
also that $\beta \in R_{m}$ for all $m \in \om$. Hence
$F(\beta) = \om$ and therefore for any $y' \in F(\beta)$ we have
that $|y - y'| \geq \ep$.
\qed

\subsection*{The case of analytic sets.}

One can ask what is the best that we can say about the set of
points of continuity of $F$ without any additional topological
assumptions for $Y$ or for $F(x)$. The following proposition gives
an upper bound for the complexity of this set. (Notice though that
we restrict ourselves to complete and separable metric spaces.)

\begin{prop}

\label{proposition the set of points of continuity is analytic}

Let $(X,p)$ and $(Y,d)$ be complete and separable metric spaces
and let $F: X \Rightarrow Y$ be a multi-valued function such that the set $F \subseteq X \times Y$ is analytic. Then the set of points of continuity of $F$ is analytic as well.
\end{prop}

\proof
Following the proof of Theorem \ref{theorem G-d when compact and
Sigma-3 when sigma-compact} for all $n \in \om$ we define the set
$A_n \subseteq X \times Y$ as follows
\[
(x,y) \in A_n \eq \inf\set{\sup\set{d(y,F(x'))}{x'\in
B_p(x,\delta)}}{\delta > 0} < \frac{1}{n+1}.
\]
Also we denote by $P$ the set of points of continuity of $F$. It
is clear that
\begin{eqnarray*}
x \in P &\eq& (\exists y \in F(x))(\forall n)[ (x,y) \in A_n]\\
&\eq& (\exists y)[ \ (x,y) \in F \ \& \ (\forall n)[(x,y) \in A_n] \ ] \ \ \ \ \ \ \ \ (1)
\end{eqnarray*}
We claim that every set $A_n$ is open. This is
again straightforward from the definitions and the fact that
\[
|d(y_1,B) - d(y_2,B)| \leq d(y_1,y_2)
\]
for all non-empty sets $B \subseteq Y$. Fix \n \ and let $(x,y) \in A_n$.
Choose $\ep_0, \delta_0 > 0$ such that
\[
\sup\set{d(y,F(x'))}{x'\in B_p(x,\delta_0)} + \ep_0 <
\frac{1}{n+1}
\]
We will show that for all $(x_1,y_1) \in X \times Y$ with $x_1 \in
B_p(x,\delta_0)$ and $y_1 \in B_d(y,\ep_0)$ we have that
$(x_1,y_1) \in A_n$. Let $x_1 \in B_p(x,\delta_0)$ and $y_1 \in
B_d(y,\ep_0)$; choose $\delta_1 > 0$ such that $B_p(x_1,\delta_1)
\subseteq B_p(x,\delta_0)$. We compute
\begin{eqnarray*}
&& \inf\set{\sup\{d(y_1,F(x'))}{x'\in B_p(x_1,\delta)\} \ / \
\delta >0} \\
&\leq& \sup\set{d(y_1,F(x'))}{x'\in B_p(x_1,\delta_1)}\\
&\leq& \sup\set{d(y_1,F(x'))}{x'\in B_p(x,\delta_0)} \\
&\leq& \sup\set{d(y,F(x'))}{x'\in B_p(x,\delta_0)} + d(y_1,y) \\
&<& \sup\set{d(y,F(x'))}{x'\in B_p(x,\delta_0)} + \ep_0 <
\frac{1}{n+1}
\end{eqnarray*}
i.e., $(x,y) \in A_n$ and the set $A_n$ is open. Now define $A = \cap A_n$; from $(1)$ above we have that
\[
x \in P \eq (\exists y)[ \ (x,y) \in F \cap A \ ],
\]
so that the set $P$ is the projection along $Y$ of the analytic set $F \cap A$. It follows that $P$ is analytic.
\qed

Now we show that if we remove just one of our assumptions about $F(x)$ or about $Y$ in Theorem \ref{theorem G-d when compact and Sigma-3 when
sigma-compact}, then it is possible that the set of
points of continuity of $F$ is not even a Borel set. Therefore
Proposition \ref{proposition the set of points of continuity is
analytic} is the best that one can say in the general case.

\begin{thm}

\hfill

\label{theorem counterexample for analytic}

\begin{enumerate}[\em(a)]

\item There is a multi-valued function $F: \ca{C} \Rightarrow \ca{N}$
such that the set $F(x)$ is closed for all $x \in \ca{C}$ and the
set of points of continuity of $F$ is analytic and not Borel. Moreover the set $F$ is a Borel subset of $\ca{C} \times \ca{N}$.

\item There is a multi-valued function $F: [0,1] \Rightarrow [0,1]$ for
which the set of the points of continuity of $F$ is analytic and not
Borel. Moreover the set $F$ is a Borel subset of $[0,1] \times [0,1]$.
\end{enumerate}

\end{thm}
\noindent
Before proving this theorem it is perhaps useful to make the following remarks. If we replace in (a) of Theorem \ref{theorem G-d when compact and Sigma-3 when sigma-compact} the condition about $F(x)$ being compact for all $x$  with ``$F(x)$ is closed for all
$x$", then from (a) of Theorem \ref{theorem counterexample for analytic} we can see that the result fails in the worst possible way. Also -in
connection with (b) of Theorem \ref{theorem G-d when compact and
Sigma-3 when sigma-compact}- we can see that if we drop the
hypothesis about $Y$ being exhaustible by compact sets but keep
the second condition ``$F(x)$ is closed for all $x$", then again
the result fails in the worst possible way.

If we replace in (a) of Theorem
\ref{theorem G-d when compact and Sigma-3 when sigma-compact} the hypothesis ``$F(x)$ is compact for all $x$", with
``$Y$ is compact" then still the result fails in the worst
possible way.

In conclusion if we want to obtain that
the set of points of continuity of a multi-valued function $F$ is
Borel, then we cannot drop the condition ``$F(x)$ is closed for
all $x$". But yet this condition alone is not sufficient in order to
derive this result as long as $Y$ is neither compact nor
exhaustible by compact sets.

\proof
We begin with (a). Let $Tr$ be the set of all (non-empty) trees on
\om, (see the Introduction). As we mentioned before the set $Tr$
can be regarded as a compact subspace of the Cantor space \ca{C}.
From Lemma \ref{lemma about extension} it is enough to construct a
multi-valued function $F: Tr \Rightarrow \ca{N}$ such that the set
of points of continuity of $F$ is not Borel, the set $F(T)$ is
closed for all $T \in Tr$ and the set $F$ is a Borel subset of $Tr \times \ca{N}$.

Denote by $IF$ the set of all ill founded trees i.e, the set of
all $T \in Tr$ for which the body $[T]$ is not empty. It is well
known (see \cite{kechris} 27.1) that the set $IF$ is an analytic subset of
$Tr$ which is not Borel.\footnote{A classical way for proving that
a given set $A \subseteq \ca{X}$ is not Borel is finding a Borel
function $\pi: Tr \to \ca{X}$ such that $IF = \pi^{-1}[A]$. If
$A$ was a Borel set then $IF$ would be Borel, a contradiction.} For
$T \in Tr$ we define the tree
\[
T^{+1} = \set{(u(0)+1, \dots, u(n-1)+1)}{u \in T, lh(u) = n}.
\]
Also we define the set $trm(T)$ as the set of all \emph{terminal}
nodes of $T$ i.e., the set of all those $u$'s in $T$ for which
there is no $w \in T$ such that $u \sqsubseteq w$ and $u \neq w$.
Define the multi-valued function $F: Tr \Rightarrow \ca{N}$ as
follows
\[
F(T) = [T^{+1}] \cup \set{u \ \hat{} \ (0,0,0,\dots)}{u \in
trm(T^{+1})}
\]
for all $T \in Tr$. Notice that if some $T \in Tr$ has no terminal
nodes then $T$ must have at least one infinite branch i.e., if
$trm(T) = \emptyset$ then $[T] \neq \emptyset$. So $F(T) \neq
\emptyset$ for all $T \in Tr$ i.e., the function $F$ is total.

First we prove that $F$ is a Borel subset of $\ca{C} \times \ca{N}$. Indeed compute
\begin{eqnarray*}
(T,\beta) \in F &\eq& \textrm{either} \ [ \ (\forall n) \ \beta(n) \geq 1 \ \& \ (\forall n) \ (\beta(0)-1,\dots,\beta(n-1)-1) \in T \ ] \\
&& \textrm{or} \ [ \ \ (\exists n)(\forall m \geq n)(\forall i < n)\\
&& \{ \ \beta(m)=0 \ \& \ \beta(i) \geq 1
\& \ (\beta(0)-1,\dots,\beta(n-1)-1) \in T
\\
&& \hspace*{1cm} \& (\forall u)[(u \in T \ \& \ lh(u)>n) \rightarrow (\exists j < n)[\beta(j) - 1 \neq u(j)] ] \ \} \  \ ].
\end{eqnarray*}
(The last line is to say that $(\beta(0)-1,\dots,\beta(n-1)-1)$ is terminal in $T$). It follows that $F$ is a $\Sigma^0_2$ set. From this and Proposition \ref{proposition the set of points of continuity is analytic} we get that the set of points of continuity of $F$ is analytic.

Now we prove that for all $T \in Tr$ the set $F(T)$ is a
closed subset of the Baire space \ca{N}. Let $(\alpha_i)_{\iin}$
be a sequence in $F(T)$ which converges to some $\alpha$. We will
prove that $\alpha \in F(T)$. If $\alpha_i \in [T^{+1}]$ for
infinitely many $i$'s then since $[T^{+1}]$ is a closed set we
have that $\alpha \in [T^{+1}]$ as well. So assume that
\[\alpha_i = u_i \ \hat{} \ (0,0,\dots) \ \ \ \textrm{for some $u_i \in
trm(T^{+1})$ for all \iin.} \ \ \ \ (1)
\]
We distinguish cases. In the first case we have that
$\sup\set{lh(u_i)}{\iin} < \infty$. Then there is a subsequence
$(u_{k_i})_{\iin}$ and some $n_0 \in \om$ such that $lh(u_{k_i}) =
n_0$ for all \iin. Hence for all \iin
\[
\alpha_{k_i}(t) = u_{k_i}(t) \ \textrm{for all $t < n_0$ and} \
\alpha_{k_i}(t) = 0 \ \textrm{for all $t \geq n_0$} \ \ \ \ (2)
\]
Since the sequence $(\alpha_{k_i}(t))_{\iin}$ is convergent for
all $t \in \om$ it follows from $(2)$ that the sequence
$(u_{k_i}(t))_{\iin}$ is also convergent for all $t < n_0$. Define
$u(t) = \lim_{\iin} u_{k_i}(t)$ for all $t < n_0$. Clearly $u$ is
a finite sequence of natural numbers with length equal to $n_0$
\footnote{If $n_0 = 0$ then $u_{k_i}$ is the empty sequence for
all $i \geq i_0$. In this case we let $u$ be the empty sequence as
well.} and there is some $i_0$ such that for all $i \geq i_0$ and
for all $t < n_0$ we have that $u_{k_i}(t) = u(t)$ i.e.,
\[
u_{k_i} = u \ \textrm{for all $i \geq i_0$.} \ \ \ \ (3)
\]
Since $u_{k_{i_0}} \in trm(T^{+1})$ we have that $u \in
trm(T^{+1})$ as well. We will show that $\alpha = u \ \hat{} \
(0,0,\dots)$. Let $n \in \om$; since $\alpha_{k_i} \to \alpha$
then there is some $i_1 \geq i_0$ such that $\alpha_{k_{i_1}}(n)
=\alpha(n)$. If $n < n_0 =lh(u)$ then -using $(2)$ and $(3)$- we have
that
\[
\alpha(n) = \alpha_{k_{i_1}}(n) = u_{k_{i_1}}(n) = u(n)
\]
and if $n \geq n_0$ then again from $(2)$
\[
\alpha(n) = \alpha_{k_{i_1}}(n) = 0.
\]
Hence $\alpha = u \ \hat{} \ (0,0,\dots)$ and since $u \in
trm(T^{+1})$ we also have that $\alpha \in F(T)$. In the second
case we have that $\sup\set{lh(u_i)}{\iin} = \infty$. Choose a
subsequence $(u_{k_i})_{\iin}$ such that $lh(u_{k_i}) \geq i$ for
all \iin. Now we claim that $\alpha \in [T^{+1}]$. Let $n \in
\om$; since $\alpha_{k_i} \to \alpha$ there is some $i_0 \in \om$
such that for all $i \geq i_0$ and for all $t < n$ we have that
$\alpha_{k_i}(t) = \alpha(t)$. For $i
> n > t$ since $lh(u_{k_i}) \geq k_i \geq i
> n > t$ from $(1)$ we have that $\alpha_{k_i}(t) = u_{k_i}(t)$.
So if we fix some $i_1 > \max\{i_0,n\}$ then for all $t < n$ we
have that $\alpha(t) = \alpha_{k_{i_1}}(t) = u_{k_{i_1}}(t)$
i.e.,
\[
(\alpha(0), \dots, \alpha(n-1)) = (u_{k_{i_1}}(0), \dots,
u_{k_{i_1}}(n-1)) \in T^{+1} \ \ \ \textrm{(because $u_{k_{i_1}}
\in T^{+1}$)}.
\]
So we have proved that $(\alpha(0), \dots, \alpha(n-1)) \in
T^{+1}$ for all \n \ i.e., $\alpha \in [T^{+1}]$.

Thus the set $F(T)$ is closed for all $T \in Tr$. Now we prove
that for $T \in Tr$
\[
\textrm{$F$ is continuous at $T$} \ \eq \ T \in IF.
\]
Since $IF$ is not Borel if we prove the above equivalence we are
done. Assume that $T$ is not in $IF$ and so $[T] = [T^{+1}] =
\emptyset$. Let any $\alpha \in F(T)$; then there is some $u$
which is terminal in $T$ such that
\[
\alpha = (u(0) + 1,\dots, u(n_0-1)+1) \ \hat{} \ (0,0,0,\dots)
\]
where $n_0 = lh(u)$. Define
\[
V = \set{\beta \in \ca{N}}{\forall i \leq n_0 \ \beta(i) =
\alpha(i)}.
\]
Clearly $V$ is a basic neighborhood of $\alpha$. Let $W$ be any
basic neighborhood of $T$ and let $u_0, \dots, u_{n-1} \in T$, $w_0,\dots w_{m-1} \in \omseq \setminus T$ be such that
\[
W = \set{T' \in Tr}{(\forall i < n)[ u_i \in T' ] \ \& \ (\forall
j < m) [ w_j \not \in T']}.
\]
We will prove that there is some $T' \in W$ such that for all
$\alpha' \in F(T')$ we have that $\alpha' \not \in V$. Since $T
\in W$ and $u \in T$ we may assume (by moving to a subset of $W$
if necessary) that $u_0 = u$. Also we may assume that for all $1
\leq i < n$ the finite sequence $u_i$ is not an initial segment of
$u_0$; (for otherwise we just remove it from the list - clearly the set $W$ is not affected).

Notice that since $T \in W$ it is not possible to have $w_j \sqsubseteq u_i$ for any $i < n$ and $j < m$.

Fix some $t_0 \in \om$ such that $u_0 \ \hat{} \ (t_0) \neq w_j$
for all $j < m$. It follows that for all $j < m$ the finite
sequence $w_j$ is not an initial segment of $u_0 \ \hat{} \
(t_0)$; for otherwise we would have that $w_j \sqsubseteq u_0$ for
some $j < m$, a contradiction. Let $T'$ be the tree which is
generated by the set $\{u_0 \ \hat{} \ (t_0), u_1, \dots,
u_{n-1}\}$ (look at the Introduction). Clearly the tree $T'$ is
finite and $u_i \in T'$ for all $i < n$. If it where $w_j \in T'$
for some $j < m$ then since $w_j$ is not an initial segment of
$u_0 \ \hat{} \ (t_0)$ there would be some $1 \leq i < n$ such
that $w_j \sqsubseteq u_i$, a contradiction. Hence $w_j \not \in
T'$ for all $j < m$. It follows that $T' \in W$.

Now let any $\alpha' \in F(T')$. Since $T'$ has empty body there
is some $u'$ which is terminal in $T'$ such that
\[
\alpha' = (u'(0) + 1,\dots, u'(n_1-1)+1) \ \hat{} \ (0,0,0,\dots)
\]
where $n_1 = lh(u')$. Since $u'$ is terminal in $T'$ it follows
that $u' \in \{u_0 \ \hat{} \ (t_0), u_1, \dots, u_{n-1}\}$. If
$u' = u_0 \ \hat{} \ (t_0) = u \ \hat{} \ (t_0)$ then $u'(n_1-1) =
t_0$ and also $lh(u') = lh(u)+1$ i.e., $n_1 = n_0 + 1$. It follows
that $\alpha'(n_0) = \alpha'(n_1-1) = u'(n_1-1) + 1 =  t_0 + 1>
0$. But on the other hand $\alpha(n_0) = 0$ and so $\alpha' \not
\in V$. If $u'= u_i$ for some $1 \leq i < n$ then since $u_i$ is
not an initial segment of $u_0 = u$ and $u_0$ is terminal in $T$ we have that the finite sequences $u$ and $u'$ are incompatible. Hence there some $k < \min\{lh(u'),
lh(u)\}$ such that $u(k) \neq u'(k)$. Since $k < \min\{lh(u'),
lh(u)\}$ we have that $\alpha'(k) = u'(k) + 1 \neq u(k) + 1 =
\alpha(k)$ i.e., there is some $k < n_0 = lh(u)$ such that
$\alpha'(k) \neq \alpha(k)$. Hence $\alpha' \not \in V$.

Now assume that $T \in IF$. We will prove that $F$ is continuous
at $T$. Let $\alpha$ be any member of $[T^{+1}]$ and let $V$ be a
basic neighborhood of $\alpha$ i.e., for some fixed $n_0$ we have
that
\[
V = \set{\beta \in \ca{N}}{\forall i \leq n_0 \ \beta(i) =
\alpha(i)}.
\]
Define
\[
W = \set{T \in Tr}{(\alpha(0) - 1, \dots, \alpha(n_0)-1) \in T}.
\]
Since $\alpha \in [T^{+1}]$ we have that $T \in W$. We will prove
that for all $T' \in W$ there is some $\alpha' \in F(T')$ such
that $\alpha' \in V$. Indeed let $T' \in W$, then $(\alpha(0) - 1,
\dots, \alpha(n_0)-1) \in T'$. There are two cases: either there
is some $u \in T'$ with $(\alpha(0) - 1,
\dots, \alpha(n_0)-1) \sqsubseteq u$ and $u$ is terminal in $T'$, or there is some $\beta \in
[T']$ such that $(\alpha(0) -1 , \dots, \alpha(n_0)-1) =
(\beta(0), \dots, \beta(n_0))$.

Assume the first case i.e., there
is some $u \in T'$ such that $(\alpha(0) - 1,
\dots, \alpha(n_0)-1) \sqsubseteq u$ and $u$ is terminal in $T'$. Clearly $lh(u) \geq n_0+1$. Take
\[
\alpha' = (u(0)+1,\dots, u(n_1-1)+1) \ \hat{} \ (0,0, \dots)
\]
where $n_1 = lh(u) \geq n_0+1$. Since $u$ is terminal in $T'$ we
have that $\alpha' \in F(T')$. Also for all $i \leq n_0 < n_1$
we have that $\alpha'(i) = u(i) + 1$; since $(\alpha(0) - 1,
\dots, \alpha(n_0)-1) \sqsubseteq u$ it follows that $\alpha(i) -
1 = u(i)$ and so $\alpha(i) = u(i) + 1 = \alpha'(i)$. Hence for all $i \leq n_0$ we have that $\alpha'(i) = \alpha(i)$ i.e., $\alpha' \in V$.

Assume the second case i.e., there is some $\beta \in
[T']$ such that $(\alpha(0) -1 , \dots, \alpha(n_0)-1) =
(\beta(0), \dots, \beta(n_0))$. Take $\alpha'(i) = \beta(i) + 1$ for all
\iin. Since $\beta \in [T']$ we have that $\alpha' \in (T')^{+1}
\subseteq F(T')$. Also for all $i \leq n_0$ we have that
$\alpha'(i) = \beta(i) + 1 = \alpha(i)$. Therefore $\alpha' \in V$
and we are done. This proves assertion $(a)$.

Assertion $(b)$ is an easy consequence of the previous. From $(a)$
and Lemma \ref{lemma about extension} it follows that there is a
multi-valued function $F: [0,1] \Rightarrow \ca{N}$ such that the
set of points of continuity of $F$ is analytic and not Borel.
Moreover the set $F$ is a Borel subset of $[0,1] \times
\ca{N}$. It is well known that there is a topological isomorphism
between $\ca{N}$ and a subset of $[0,1]$; we give a description of
its construction. By induction one constructs a family $(I_u)_{u
\in \omseq}$ of closed intervals such that $I_{\empt} = [0,1]$;
length$(I_u) \leq \frac{1}{2^{lh(u)}}$ for all $u$; $I_{u \ \hat{}
\ (n)} \subseteq I_u$ for all $u$, $n$; and $I_u \cap I_w =
\emptyset$ for all incompatible $u$, $w$. It is not hard to see
that the function $\pi: \ca{N} \to [0,1]$ defined by
$\pi(\alpha)=$ \emph{the unique $x \in [0,1]$ such that $x \in
\bigcap\limits_{n} I_{(\alpha(0),\dots, \alpha(n))}$} is a
topological isomorphism between $\ca{N}$ and a $G_\delta$ subset
of $[0,1]$.

Since there is a topological isomorphism between $\ca{N}$ and a
subset of $[0,1]$ from the first assertion of Lemma \ref{lemma
about composition} there is a multi-valued function $G: [0,1]
\Rightarrow [0,1]$ for which the set of points of continuity is
analytic and not Borel. Moreover from the third assertion of Lemma
\ref{lemma about composition} the set $G$ is a Borel subset of
$[0,1] \times [0,1]$.\footnote{Notice (in connection with Theorem
\ref{theorem G-d when compact and Sigma-3 when sigma-compact})
that $\pi[\ca{N}]$ is not a closed subset of $[0,1]$ -for
otherwise the Baire space would be compact. So the conclusion of
the second assertion of Lemma \ref{lemma about composition} does
not apply here.}
\qed

\textbf{Question 1.}
\begin{enumerate}
\item In all theorems and examples we have seen about multi-valued functions, the set of points of continuity is either a $\Sigma^0_3$ set or not even a Borel set. It would be interesting to see if one can construct multi-valued functions, whose set of points of continuity lies somewhere between the $\Sigma^0_3$-level and the level of analytic sets. More specifically: does there exist for every $n > 3$ a multi-valued function, whose set of points of continuity is a $\Sigma^0_n$ (or $\Pi^0_n$) and not a $\Pi^0_n$ ($\Sigma^0_n$ respectively) set?

    A positive answer to this question might be given by investigating the next question.

\item Suppose that we are given a
multi-valued function $F: X \Rightarrow Y$ for which we have that
the set $F(x)$ is closed for all $x$ and $Y$ is separable. As we
have proved before in case where $Y$ is exhaustible by compact
sets the set of points of continuity of $F$ is $\Sigma^0_3$ and in
case where $Y = \ca{N}$ it is possible that the latter set is not
even Borel. In fact one can see that the latter is true not just
for $Y = \ca{N}$ but also in case where \ca{N} is topologically
isomorphic with a closed subset of $Y$; (see the second assertion
of Lemma \ref{lemma about composition} and the proof of the second
assertion of Theorem \ref{theorem counterexample for analytic}).
The question is what happens when $Y$ falls in neither of the
previous cases i.e., $Y$ is neither exhaustible by compact sets
nor it contains \ca{N} as a closed subset. An interesting class of
such examples is the class of infinite dimensional separable
\emph{Banach spaces} i.e., (infinite dimensional) linear normed
spaces which are complete and separable under that norm. Any such
space is not exhaustible by compact sets and it does not contain
\ca{N} as a closed subset. Therefore the theorems of this article
provide no information in this case. It would be interesting to
find the best upper bound for the complexity of the set of points
of continuity of $F$ when $Y$ is an infinite dimensional separable
Banach space and the set $F(x)$ is closed for all $x$.

\end{enumerate}

\section{Strong Continuity.}

\label{section strong continuity}

We continue with some results regarding the set of
points of \emph{strong} continuity of a multi-valued function $F$.
In particular we will prove the corresponding of Theorem
\ref{theorem G-d when compact and Sigma-3 when sigma-compact} and
Proposition \ref{proposition the set of points of continuity is
analytic}. Examples which show that these results are optimum is a subject which is still under investigation. Before we proceed the author would like to express his special thanks to the anonymous referee who has provided him with very interesting remarks about strong continuity. These remarks include the observation that $F$ is strongly continuous at $x_0$ exactly when the multi-valued function $x \mapsto \overline{F(x)}$ (where $\overline{F(x)}$ is the topological closure of $F(x)$) is strongly continuous at $x_0$, which is continuity at $x_0$ with respect to the lower Fell topology. These remarks have motivated Proposition \ref{proposition about Fell topology}.

As we mentioned in the beginning,
Theorem \ref{theorem G-d when compact and Sigma-3 when
sigma-compact} does not require any additional hypothesis about
$F$ as a subset of $X \times Y$. However the following remark
suggests that this is not the case for strong continuity.

\begin{rem}

\label{remark about strong continuity}

\normalfont

Let $A$ be a dense subset of $[0,1]$. We define the multi-valued function $F: [0,1] \Rightarrow \{0,1\}$ as follows
\[
F(x) = \{0\}, \ \ \textrm{if $x \in A$ and} \ F(x) = \{0,1\} \ \ \textrm{if $x \not \in A$},
\]
for all $x \in [0,1]$. We claim that the set of points of strong continuity of $F$ is exactly the set $A$. Let $x \in A$, $y \in F(x)$ and $\ep > 0$. Take $\delta = 1 > 0$ and let $x' \in (x-\delta,x+\delta)$. We have that $y = 0$ and also since $0 \in F(x')$ we can take $y' = 0$; so $|y-y'| = 0 < \ep$. Now let $x \not \in A$. We take $y = 1 \in F(x)$ and $\ep = \frac{1}{2}$. Let any $\delta >0$. Since $A$ is a dense subset of $[0,1]$ there is some $x' \in A$ such that $x' \in (x-\delta,x+\delta)$. Clearly for all $y' \in F(x')$ we have that $y' = 0$ and so $|y - y'| = 1 > \ep$.\qed
\end{rem}

Since there are dense subsets of $[0,1]$ which are way above the level of analytic sets from Remark \ref{remark about strong continuity} we can see that there is no hope to obtain the corresponding of Theorem \ref{theorem G-d when compact and Sigma-3 when sigma-compact} without any additional assumptions about the complexity of the set $F$. Notice also that those assumptions about the set $F$ have to be at least as strong as the result that we want to derive. For example it is well known that there is a dense $\Pi^0_3$ set $A \subseteq [0,1]$ which is not $\Sigma^0_3$;  hence by taking the multi-valued function $F$ of Remark \ref{remark about strong continuity} with respect to that set $A$ we can see that $F$ is $\Delta^0_4$ as a subset of $[0,1] \times [0,1]$ and that the set of points of strong continuity of $F$ (i.e., the set $A$) is not $\Sigma^0_3$. In other words if we want to result to a $
\Sigma^0_3$ set we need to assume that $F$ does not go above the third level of the Borel hierarchy. The following may be regarded as the corresponding strong-continuity analogue of Theorem \ref{theorem G-d when compact and Sigma-3 when sigma-compact}.

\begin{thm}

\label{theorem strong continuity G-d when compact and Sigma-3 when sigma-compact}

Let $(X,p)$ and $(Y,d)$ be metric spaces with $(Y,d)$ being
separable and let $F: X \Rightarrow Y$ be a multi-valued function such that $F$ is a $\Sigma^0_2$ subset of $X \times Y$.

\begin{enumerate}[\em(a)]

\item If $Y$ is compact and the set $F(x)$ is closed for all $x \in X$ then the set of points of strong continuity of $F$ is $\Pi^0_2$.

\item If $Y$ is exhaustible by compact sets
and the set $F(x)$ is closed for all $x \in X$, then the set of
points of strong continuity of $F$ is $\Sigma^0_3$.
\end{enumerate}

\end{thm}

\proof
Following the proof of Theorem \ref{theorem G-d when compact and
Sigma-3 when sigma-compact} we denote by $P$ the set of $x$'s for
which $F$ is strongly continuous at $x$ and also we fix a sequence
$(y_s)_{s \in \om}$ in $Y$ for which the set $\set{y_s}{\s}$ is
dense in $Y$. Let us begin with (a); we claim that
\begin{eqnarray*}
x \in P &\eq& (\forall n)(\forall s)[ \ \ d(y_s,F(x)) > \frac{1}{3(n+1)} \\
&& \hspace*{1.2cm} \textrm{or} \
\inf\set{\sup\set{d(y_s,F(x'))}{x'\in B_p(x,\delta)}}{\delta > 0}
< \frac{1}{n+1} \ ] \ \ \ \ \ \ \ (*)
\end{eqnarray*}
The proof is similar to the corresponding equivalence in the proof
of Theorem \ref{theorem G-d when compact and Sigma-3 when
sigma-compact}. We give a brief description. For the
left-to-right-hand direction notice that if $d(y_s,F(x)) \leq
\frac{1}{3(n+1)} < \frac{1}{2(n+1)}$ then there is some $y \in
F(x)$ such that $d(y_s,y) < \frac{1}{2(n+1)}$; then proceed
exactly as before. For the converse direction let $y \in F(x)$ and
choose a sequence $(y_{s_n})_{\n}$ such that $d(y_{s_n},y) <
\frac{1}{3(n+1)}$ for all \n. From the hypothesis we have that
\[
\inf\set{\sup\set{d(y_{s_n},F(x'))}{x'\in B_p(x,\delta)}}{\delta
> 0} < \frac{1}{n+1}
\]
for all \n. Following the same steps as before we get the result.
Now define the sets $Q_{n,s}$ and $R_{s,n}$ as follows:
\begin{eqnarray*}
x \in Q_{n,s} &\eq& \inf\set{\sup\set{d(y_s,F(x'))}{x'\in B_p(x,\delta)}}{\delta > 0} < \frac{1}{n+1},\\
x \in R_{s,n} &\eq& d(y_s,F(x)) \leq \frac{1}{3(n+1)}.
\end{eqnarray*}
From $(*)$ we have that $P = \bigcap\limits_{n,s} [(X \setminus
R_{n,s}) \cup Q_{n,s}]$. We have already proved that the set
$Q_{n,s}$ is open for all $n,s$, so it is enough to prove that
$R_{n,s}$ is $\Sigma^0_2$ for all $n,s$. Since $F$ is $\Sigma^0_2$
we may write $F = \bigcup\limits_{j \in \om}F_j$ where $F_j$ is a
closed subset of $X \times Y$ for all $j \in \om$. We adopt the
notation $F_j(x)$ for the $x$-section of $F_j$. Using the
compactness of $F(x)$ we obtain that
\[
d(y_s,F(x)) \leq \frac{1}{3(n+1)} \eq (\exists j) d(y_s,F_j(x))
\leq \frac{1}{3(n+1)}.
\]
Define $x \in R_{n,s,j} \eq d(y_s,F_j(x)) \leq \frac{1}{3(n+1)}$.
We will prove that $R_{n,s,j}$ is closed for all $n,s,j$. Indeed
let $(x_i)_{\iin} \subseteq R_{n,s,j}$ be such that $x_i
\stackrel{\iin}{\longrightarrow}x$. Since $F_j(x_i)$ is compact
there is some $z_i \in F_j(x_i)$ such that $d(y_s,z_i) =
d(y_s,F_j(x_i)) \leq \frac{1}{3(n+1)}$ for all \iin. From the
compactness of $Y$ there is a subsequence $(z_{k_i})_{\iin}$ and
some $z \in Y$ such that $z_{k_i}
\stackrel{\iin}{\longrightarrow}z$. Since $(x_{k_i},z_{k_i}) \in
F$ for all \iin \ and $F$ is closed we have that $(x,z) \in F$.
Hence $z \in F(x)$ and so $d(y_s,F(x))\leq d(y_s,z) =
\lim_{\iin}d(y_s,z_{k_i}) \leq \frac{1}{3(n+1)}$ i.e., $x \in
R_{n,s,j}$.

For $(b)$, as previously we write $Y = \bigcup\limits_{m \in \om}
K_m$ where $(K_m)_{m \in \om}$ is a sequence of compact subsets of
$Y$ such that $K_m \subseteq K_{m+1}^{\circ}$ for all $m \in \om$.
We claim the following modified version of $(*)$:
\begin{eqnarray*}
x \in P &\eq& (\exists m)(\forall n)(\forall s) [ \ d(y_s,F(x) \cap K_m) > \frac{1}{3(n+1)} \\
&& \hspace*{0.3cm} \textrm{or} \ \inf\set{\sup\set{d(y_s,F(x')\cap
K_m)}{x'\in B_p(x,\delta)}}{\delta
> 0} < \frac{1}{n+1} \ ] \ \ \ \ \ \ (\dagger)
\end{eqnarray*}
This follows from the proof of Theorem \ref{theorem G-d when
compact and Sigma-3 when sigma-compact} with the modifications we
described before. The relation $R^m_{n,s,j}$ defined by
\[
x \in R^m_{n,s,j} \eq d(y_s,F_j(x) \cap K_m) \leq \frac{1}{3(n+1)}
\]
is closed; the proof is as above. (Notice that if $x \in
R^m_{n,s,j}$ we have that $F_j(x) \cap K_m \neq \emptyset$.)
\qed

While Remark \ref{remark about strong continuity} makes it clear that the notion of strong continuity is non-metrizable in the sense of Corollary \ref{corollary continuity is non-metrizable}, the notions of strong continuity and metrizability are not entirely unrelated. First we put down the necessary definitions. Suppose that $Y$ is a topological space. We denote the family of all closed subsets of $Y$ by $\ca{F}(Y)$. The \emph{lower Fell topology} on $\ca{F}(Y)$ is the topology generated by all sets of the form $\ca{A}_U= \set{C \in \ca{F}(Y)}{C \cap U \neq \emptyset}$, where $U$ is an open subset of $Y$. We now consider the least $\sigma$-algebra $S$ on $\ca{F}(Y)$ containing all sets of the form $\ca{A}_U$. The pair $(\ca{F}(Y),S)$ is the \emph{Effros Borel space} of $\ca{F}(Y)$. A well-known theorem states that if $Y$ is separable and metrizable by a complete distance function, (i.e., $Y$ is a \emph{Polish space}) then there is a topology $\ca{T}$ on $\ca{F}(Y)$ such that: (a) the space $(\ca{F}(Y),\ca{T})$ is a Polish space and (b) the $\ca{T}$-Borel subsets of $\ca{F}(Y)$ are exactly the members of $S$, c.f. \cite{kechris} Section 12.C. The \emph{Fell topology} on $Y$ is the topology which has as basis the family of all sets of the form $\ca{W} \equiv \ca{W}(K,U_1,\dots,U_n) = \set{C \in \ca{F}(Y)}{C \cap K = \emptyset\ \& \ (\forall i \leq n)[C \cap U_i \neq \emptyset]}$, where $K$ is a compact subset of $Y$ and $U_1, \dots, U_n$ are open subsets of $Y$. By choosing $K$ as the empty set one can see that the lower Fell topology is contained in the Fell topology. If $Y$ is a locally compact Polish space then the Fell topology is compact metrizable and its Borel space is exactly the Effros Borel space, c.f. \cite{kechris} Section 12.C. Suppose that $X$ is a metric space, $Y$ is a Polish space and $F$ is a multi-valued function from $X$ to $Y$. We consider the multi-valued function $\overline{F}: X \Rightarrow Y: \overline{F}(x)=\textrm{the topological closure of $F(x)$}$. We say that $F$ is \emph{pointwise-closed} exactly when $F(x)$ is closed for all $x$, which is of course the same as saying that $F = \overline{F}$.

\begin{prop}

\label{proposition about Fell topology}

\normalfont

Consider a multi-valued function $F: (X,p) \Rightarrow (Y,d)$, with $(Y,d)$ complete and separable. The following hold.
\begin{enumerate}[(1)]

\item The multi-valued function $F$ is strongly continuous at $x \in X$ exactly when the multi-valued function $\overline{F}$ is strongly continuous at $x$, \footnote{Unfortunately this observation does not seem to allow us to drop the hypothesis ``$F(x)$ is closed" in Theorem \ref{theorem strong continuity G-d when compact and Sigma-3 when sigma-compact}.} which is exactly when the function $\overline{F}: X \to \ca{F}(Y)$ is continuous at $x$ with respect to the lower Fell topology on $\ca{F}(Y)$. \emph{In particular if $\overline{F}: X \to \ca{F}(Y)$ is continuous at $x$ with respect to the Fell topology then $F$ is strongly continuous at $x$.}
\item Consider a Polish topology $\ca{T}$ on $\ca{F}(Y)$ such that the Borel space with respect to \ca{T} is the Effros Borel space. If $F$ is strongly continuous at every $x \in X$, then the function $\overline{F}: X \to (\ca{F}(Y),\ca{T})$ is Borel measurable. (We are thinking of $X$ with its Borel structure.) \emph{Thus strong continuity is reduced to Borel measurability for pointwise-closed multi-valued functions with range a Polish space.}

\item If $F$ is a closed subset of $X \times Y$ and it is strongly continuous at $x \in X$ then the function $F: X \to \ca{F}(Y)$ is continuous at $x$ with respect to the Fell topology.

\item Suppose $F$ is a closed subset of $X \times Y$. Then the multi-valued function $F: X \Rightarrow Y$ is strongly continuous at $x \in X$ exactly when the function $F: X \to \ca{F}(Y)$ is continuous at $x$ with respect to the Fell topology.

    \emph{It follows that in the case of multi-valued functions with closed graph and range a locally compact Polish space, the notion of strong continuity is metrizable.} In particular in any example which shows that the result of Theorem \ref{theorem strong continuity G-d when compact and Sigma-3 when sigma-compact}-(b) is optimal, the set $F$ must be $\Sigma^0_2$ and not $\Pi^0_1$.
\end{enumerate}
\end{prop}

\proof
\hfill
\begin{enumerate}[(1)]
\item The first assertion is clear. Assume now that $\overline{F}$ is strongly continuous at $x$; let $U$ be open in $Y$ and such that $\overline{F}(x) \in \ca{A}_{U}$ i.e., $\overline{F}(x) \cap U \neq \emptyset$. Pick some $y \in \overline{F}(x)$ with $y \in U$ and $\ep > 0$ with $B_d(y,\ep) \subseteq U$. Since $\overline{F}$ is strongly continuous at $x$, there exists some $\delta > 0$ such that for all $x' \in B_p(x,\delta)$ there exists some $y' \in \overline{F}(x')$ such that $y' \in B_d(y,\ep)$. Hence $\overline{F}(x') \cap U \neq \emptyset$ for all $x' \in B_p(x,\delta)$. Conversely assume that $\overline{F}$ is continuous at $x$ with respect to the lower Fell topology; let $y \in \overline{F}(x)$ and $\ep > 0$. We consider the open neighborhood $\set{C \in \ca{F}(Y)}{C \cap B_d(y,\ep) \neq \emptyset}$ of $\overline{F}(x)$. There exists some $\delta > 0$ such that for all $x' \in B_p(x,\delta)$ we have that $\overline{F}(x') \cap B_d(y,\ep) \neq \emptyset$. This means that there exists some $\delta > 0$ such that for all $x' \in B_p(x,\delta)$ there exists some $y' \in \overline{F}(x')$ with $d(y,y') < \ep$.

    The last assertion of 1. follows from the fact that the lower Fell topology is contained in the Fell topology.

\item This is an immediate consequence of 1. and the comments about the Effros Borel space proceeding Proposition \ref{proposition about Fell topology}.

\item This is similar to the proof of 1. Suppose that 
\[
\ca{W} = \set{C \in \ca{F}(Y)}{C \cap K = \emptyset\ \& \ (\forall i \leq n)[C \cap U_i \neq \emptyset]}
\]
  is a basic neighborhood of $F(x)$ with respect to the Fell topology. Choose $y_i \in F(x) \cap U_i$ and $\ep_i > 0$ with $B_d(y_i,\ep_i) \subseteq U_i$ for all $i=1,\dots,n$. From our hypothesis there is some $\delta_0 >0$ such that for all $x' \in B_p(x,\delta_0)$ we have that $F(x') \cap B_d(y_i,\ep_i) \neq \emptyset$ for all $i=1,\dots,n$.

  Now we claim that there is some $0 < \delta < \delta_0$ such that for all $x' \in B_p(x,\delta)$ we have that $F(x') \cap K = \emptyset$. Indeed if this is not the case, then there are sequences $(x_n)_{\n}$ and $(z_n)_{\n}$ such that $x_n \to x$ and $z_n \in F(x_n) \cap K$ for all \n. From the compactness of $K$ there is a subsequence $(z_{k_n})_{\n}$ which converges to some $z \in K$. Therefore $(x_{k_n},z_{k_n}) \to (x,z)$. Since $(x_{k_n},z_{k_n}) \in  F$ for all \n \ and $F$ is closed it follows that $z \in F(x)$. Thus $F(x) \cap K \neq \emptyset$ contradicting that $\ca{W}$ is a neighborhood of $F(x)$.

  It follows that there is some $\delta > 0$ such that for all $x' \in B_p(x,\delta)$ we have that $F(x') \in \ca{W}$.

\item The first assertion is clear from 1. and 3. The second assertion follows from the fact that when $Y$ is Polish and locally compact the Fell topology is compact metrizable. The assertion about Theorem \ref{theorem strong continuity G-d when compact and Sigma-3 when sigma-compact} follows from the fact that every space which is exhaustible by compact sets is locally compact.
\qed
\end{enumerate}

\noindent
We continue with the corresponding of Proposition \ref{proposition the set of points of continuity is analytic}.

\begin{prop}

\label{proposition strong continuity the set of points of continuity is co-analytic}

Let $(X,p)$ and $(Y,d)$ be complete and separable metric spaces
and let $F: X \Rightarrow Y$ be a multi-valued function such that the set $F \subseteq X \times Y$ is analytic. Then the set of points of strong continuity of $F$ is co-analytic.
\end{prop}

\proof
As in the proof of Proposition \ref{proposition the set of points of continuity is analytic} for all $n \in \om$ we define the set
$A_n \subseteq X \times Y$ as follows
\[
(x,y) \in A_n \eq \inf\set{\sup\set{d(y,F(x'))}{x'\in
B_p(x,\delta)}}{\delta > 0} < \frac{1}{n+1}.
\]
As before the set $A_n$ is open for all \n. Let $P_s$ be the set of points of strong continuity of $F$. It follows that
\begin{eqnarray*}
x \in P_s &\eq& (\forall y \in F(x))(\forall n)[ (x,y) \in A_n]\\
&\eq& (\forall y)[ \ (x,y) \not \in F \ \textrm{or} \ (\forall n)[(x,y) \in A_n] \ ].
\end{eqnarray*}
So
\begin{eqnarray*}
x \not \in P_s \eq (\exists y)[ \ (x,y) \in F \ \& \ (\exists n)[(x,y) \not \in A_n] \ ].
\end{eqnarray*}
Since $F$ is analytic and since every set $A_n$ is open we obtain that the complement of $P_s$ is analytic i.e., the set $P_s$ is co-analytic.
\qed

We conclude this article with some remarks which concern all previous results. The author would like to thank the anonymous referee for raising the questions stated below.

\begin{rem}

\normalfont

\label{remark about lightface}

All results above are in the context of classical descriptive set theory. One could ask whether the corresponding results are also true in the context of \emph{effective} descriptive set theory. In the latter context one deals with the notion of a recursive function $f: \om^k \to \om^m$ and of a recursive subset of $\om^k$. We assume that our given metric space $(X,d)$ is complete, separable and that there is a countable dense sequence \set{r_i}{i \in \om} such that the relations $d(r_i,r_j) < q$, $d(r_i,r_j) \leq q$ for $i,j \in \om$ and $q \in \mathbb{Q}^+$, are recursive. (An example of such space is $\mathbb{R}$ with $\set{r_i}{i \in \om}=\mathbb{Q}$.) One takes then the family \set{N(X,s)}{\s} of all open balls with centers from the set \set{r_i}{i \in \om} and rational radii and defines the class of \emph{semirecursive} sets or ``effectively open" sets as the sets which are recursive unions of sets of the form $N(X,s)$. The analogous notions go through the whole hierarchy of Borel and analytical sets i.e., one constructs the family of effectively closed, effectively $G_\delta$, effectively analytic sets and so on. The latter classes of sets are also called \emph{lightface} classes. The usual inclusion properties hold also for the lightface classes. For example every effectively closed set is effectively $G_\delta$. We should point out that there are only countably many subsets of a fixed space $X$ which belong to a specific lightface class. Also all singletons $\{q\}$ with $q \in \mathbb{Q}$ belong to every one of the lightface classes mentioned above except from the one of semi-recursive sets. The reader can refer to \cite{yiannis} for a detailed exposition of this theory. One natural question which arises is if the results which are presented in this article hold in the context of effective descriptive set theory. For example: if $F: \mathbb{R} \Rightarrow \mathbb{R}$ is a bounded multi-valued function such that the set $F(x)$ is effectively closed, is it true that the set of points of continuity of $F$ is effectively $G_\delta$? As the next proposition shows the answer to this question is negative even if $F$ is a single-valued function.
\end{rem}

Let us say that a family of sets \G \ is \emph{closed under negation} if whenever $A \subseteq X$ is in \G \ then $X \setminus A$ is in \G \ as well.

\begin{prop}

\label{proposition counterexample for effective}

Suppose that \G \ is a class of sets which is closed under negation and the family \newline \set{A\subseteq \mathbb{R}}{A \in \G} is countable. Then there is a function $f: \mathbb{R} \to \set{\frac{1}{n+1}}{\n} \cup \{0\}$ such that the set of points of continuity of $f$ is not a member of \G. In particular (by choosing \G \ as the lightface $\Delta^1_2$ class) there is a function $f: \mathbb{R} \to [0,1]$ such that the singleton $\{f(x)\}$ is effectively closed for all $x \in \mathbb{R}$ but the set of points of continuity of $f$ is neither effectively analytic nor effectively co-analytic.
\end{prop}

\proof
Since \G \ restricted on the subsets of $\mathbb{R}$ is countable there is some $A \subseteq \mathbb{R}$ which is infinite countable set and not a member of \G. Write $A = \set{x_n}{\n}$ with $x_n \neq x_m$ for $n \neq m$ and define $f: \mathbb{R} \to \mathbb{R}$ as follows: $f(x) = 0$ if $x \not \in A$ and $f(x) = \frac{1}{n+1}$ if $x = x_n$.

We claim that the set of points of continuity of $f$ is $\mathbb{R} \setminus A$. Since $A$ is not a member of $\Gamma$ and the class \G \ is closed under negation we also have that the set $\mathbb{R} \setminus A$ is not a member of \G. So if we prove our claim we are done.

Assume that $x \in A$ and $x = x_n$. Choose $\ep = \frac{1}{2(n+1)} > 0$ and let any $\delta > 0$. The interval $(x-\delta,x+\delta)$ is uncountable so it contains some $y$ which is not a member of $A$. It follows that $f(y) = 0$ and so $|f(y) - f(x)| = f(x) = \frac{1}{n+1} > \ep$. This shows that no point of $A$ is a point of continuity of $f$. Now assume that $x \not \in A$ and let arbitrary $\ep > 0$. Choose \n \ such that $\frac{1}{n+1} < \ep$ and define
\[
\delta = \min\set{|x-x_i|}{i=0,\dots,n}.
\]
Since $x \not \in A$ we have that $\delta > 0$. Now let any $y \in (x-\delta,x+\delta)$; if $y \in A$ then $y = x_m$ for some $m > n$, so $|f(y) - f(x)| = f(y) = \frac{1}{m+1} < \frac{1}{n+1} < \ep$. If $y \not \in A$ then clearly $|f(y) - f(x)| = 0 <\ep$. Therefore the function $f$ is continuous at $x$.
\qed

\textbf{Question 2.} In case we take $\G$ to be the lightface $\Delta^0_n$ class for some small \n, it would be interesting to see whether one can construct a function $f$ which satisfies the first conclusion of the previous proposition and has the additional property that the graph of $f$ belongs to \G.

\bibliographystyle{eptcs}

\end{document}